\documentclass[11pt,reqno]{amsproc}

\title[On the inviscid limit of the NSE]{On the inviscid limit of the Navier-Stokes equations}

\author{Peter Constantin}
\author{Igor Kukavica}
\author{Vlad Vicol}

\address{Department of Mathematics, Princeton University, Princeton, NJ 08544}
\email{const@math.princeton.edu}

\address{Department of Mathematics, University of Southern California, Los Angeles, CA 90089}
\email{kukavica@usc.edu}

\address{Department of Mathematics, Princeton University, Princeton, NJ 08544}
\email{vvicol@math.princeton.edu}

\usepackage[margin=1in]{geometry}
\usepackage{amsmath, amsthm, amssymb}
\usepackage{color}
\usepackage{times}
\theoremstyle{plain}
\newtheorem{theorem}{Theorem}[section]

\theoremstyle{definition}
\newtheorem{remark}[theorem]{Remark}
\newtheorem{example}[theorem]{Example}

\numberwithin{equation}{section}

\renewcommand{\tilde}{\widetilde}

\def\phi{\varphi}

\def\RR{{\mathbb R}}
\def\HH{{\mathbb H}}
\def\OO{{\mathcal O}}

\def\uns{{u}}
\def\ue{{\bar{u}}}
\def\pns{{p}}
\def\pe{{\bar{p}}}
\def\uu{{U}}

\def\curl{\mathop{\rm curl} \nolimits}
\def\ddiv{\mathop{\rm div} \nolimits}

\chardef\coloryes=0
\chardef\isitdraft=0
\ifnum\isitdraft=1
   \usepackage{fancyhdr}
   \input macros.tex
   
\else

\fi

\def\comma{ {\rm ,\qquad{}} }            %comma in a formula
\def\fractext#1#2{{#1}/{#2}}
\def\dist{\mathop{\rm dist}\nolimits}    %distance
\def\indeq{\qquad}
\begin{document}

%%%%%%%%%%%%%%%%%%%%%%%%% THE ABSTRACT %%%%%%%%%%%%%%%%%%%%%%%%%%%%%%%%%%%

\begin{abstract}
We consider the convergence in the $L^2$ norm, uniformly in time, of the Navier-Stokes equations with Dirichlet boundary conditions to the Euler equations with slip boundary conditions. We prove that if the Oleinik conditions of no back-flow in the trace of the Euler flow, and of a lower bound for the Navier-Stokes vorticity is assumed in a Kato-like boundary layer, then the inviscid limit holds. 
\hfill \today.
\end{abstract}

%%%%%%%%%%%%%%%%%%%%%%%% Classification and Keywords %%%%%%%%%%%%%%%%%%%%



\maketitle

%%%%%%%%%%%%%%%%%%%%%%%%%% The Main Part %%%%%%%%%%%%%%%%%%%%%%%%%%%%%%%%%%
\section{Introduction}

We consider the two-dimensional Navier-Stokes equation \eqref{eq:NS} for the velocity field $\uns = (\uns_1,\uns_2)$ and pressure scalar $\pns$, and the two-dimensional Euler equation \eqref{eq:E} for the velocity field $\ue = (\ue_1,\ue_2)$ and scalar pressure $\pe$
 \begin{align}
  &\partial_t \uns - \nu \Delta \uns + \uns \cdot \nabla \uns + \nabla \pns = 0 \label{eq:NS}\\
  &\partial_t \ue + \ue \cdot \nabla \ue + \nabla \pe = 0 \label{eq:E}
 \end{align}
in the half plane $\HH = \{x =(x_1,x_2) \in \RR^2 \colon x_2 > 0 \}$ with Dirichlet and slip boundary conditions
 \begin{align}
  &\uns|_{\partial \HH} = 0
   \label{EQ22}
  \\&
  \ue_2|_{\partial \HH} = 0
   \label{EQ23}
 \end{align}
on the Navier-Stokes and Euler solutions respectively. 
The choice of domain being the half-plane $\HH$ is made here for simplicity of the presentation. Indeed, as discussed in Section~\ref{sec:curved} below, the results in this paper also hold if the equations are posed in a bounded domain $\Omega$ with smooth boundary. 

The initial conditions for the Euler and Navier-Stokes equations are taken to be the same,
$
\uns_0 = \ue_0. 
$
We shall also denote the Navier-Stokes vorticity as
\begin{align*}
 \omega = \partial_1 \uns_2 - \partial_2 \uns_1,
\end{align*}
and by
\begin{align*} 
 \uu=\ue_1|_{\partial \HH} 
\end{align*}
the trace of the tangential component of the Euler flow.  

Before we describe the results, we comment on scaling. We choose units of length and units of time associated to this Euler trace so that in the new variables the Euler solution $\ue$ becomes ${\mathcal O}(1)$. The integral scale ${\mathcal L}$ is given by ${\mathcal L} = \| \uu\|_{L^\infty_t L^2_x}^2 \| \uu\|_{L^\infty_{t,x}}^{-2}$ and the time scale ${\mathcal T}$ is chosen to be ${\mathcal T} = {\mathcal L} \| \uu\|_{L^\infty_{t,x}}^{-1}$. Using ${\mathcal L}$ and ${\mathcal T}$ we non-dimensionalize the Euler and Navier-Stokes equations, but for notational convenience we still refer to the resulting Reynolds number $Re = {\mathcal L}^2  {\mathcal T}^{-1} \nu^{-1}$ as $\nu^{-1}$. For the remainder of the paper the this rescaling is implicitly used, and all quantities involved are dimensionless. 

Since $\uu \neq 0$ in general, there is a mismatch between the Navier-Stokes and Euler boundary conditions leading to the phenomenon of boundary layer separation. Establishing whether
 \begin{align}
 \| \uns - \ue \|_{L^\infty(0,T;L^2(\HH))}^2 \to 0 \label{eq:INVISCID:LIMIT}
 \end{align} 
holds in the inviscid limit $\nu \to 0$ is an outstanding physically important problem in fluid dynamics. Here $T >0$ is a fixed $\nu$-independent time. There is a vast literature on the subject of inviscid limits. We refer the reader for instance to~\cite{ConstantinWu95,ConstantinWu96} for the case of vortex patches, and to~\cite{Kato84b,TemamWang97b,Masmoudi98,SammartinoCaflisch98b,OleinikSamokhin99,Wang01,Kelliher07,Masmoudi07,Kelliher08,LopesMazzucatoLopes08,
MazzucatoTaylor08,Kelliher09,Maekawa13,Maekawa14,GrenierGuoNguyen14} and references therein, for inviscid limit results in the case of Dirichlet boundary conditions. 

Going back at least to the work of Prandtl~\cite{Prandtl1904}, based on matched asymptotic expansions, one may formally argue that as $\nu \to 0$ we have
 \begin{align} 
 \uns(x_1,x_2,t) \approx \ue(x_1,x_2,t) {\bf 1}_{\{x_2>\sqrt{\nu}\}} + u_{P}(x_1,x_2/\sqrt{\nu},t) {\bf 1}_{\{x_2<\sqrt{\nu}\}} + \OO(\sqrt{\nu}) 
 \label{eq:expansion}
 \end{align}
where $u_P$ is the solution of the Prandtl boundary layer equations. We refer the reader to~\cite{Oleinik66,EEngquist97,OleinikSamokhin99,SammartinoCaflisch98a,CaflischSammartino00,E00,CannoneLombardoSammartino01,HongHunter03,Grenier00,XinZhang04,GarganoSammartinoSciacca09,GerardVaretDormy10,GuoNguyen10,GerardVaretNguyen12,MasmoudiWong12a,AlexandreWangXuYang12,KukavicaVicol13a,GerardVaretMasmoudi13,KukavicaMasmoudiVicolWong14}
for results regarding the Prandtl boundary layer equations.

The Prandtl solution is believed to describe the creation and the evolution of vorticity in a boundary layer of thickness $\sqrt{\nu}$, which makes the problem of establishing the inviscid limit \eqref{eq:INVISCID:LIMIT} intimately related to the question of well-posedness for the Prandtl equations. 

We emphasize however that up to our knowledge there is currently no abstract result which states that if the Prandtl equations are well-posed, then the inviscid limit of Navier-Stokes to Euler holds in $L^\infty(0,T;L^2)$. This is the main motivation for our paper.

So far the well-posedness of the Prandtl equation has been established in the following settings:
\begin{itemize}
\item[(a)] There is no back-flow in the initial velocity field, i.e., $U_0>0$, and the initial vorticity is bounded from below by a strictly positive constant $\omega_0 \geq \sigma >0$. This result goes back to Oleinik~\cite{Oleinik66}, and we refer to~\cite{MasmoudiWong12a} for an elegant Sobolev energy-based proof. 
\item[(b)] The initial velocity is real-analytic with respect to both the normal and tangential variables~\cite{SammartinoCaflisch98a}. 
\item[(c)] The initial velocity is real-analytic with respect to only the tangential variable~\cite{CannoneLombardoSammartino01,KukavicaVicol13a}.
\item[(d)] The initial vorticity has a single curve of non-degenerate critical points, and it lies in the Gevrey-class $7/4$ with respect to the tangential variable~\cite{GerardVaretMasmoudi13}.
\item[(e)] The initial data is of finite Sobolev smoothness, the vorticity is positive on an open strip $(x,y) \in I \times [0,\infty)$,   is negative for $(x,y) \in I^C\times[0,\infty)$, and the vorticity is real-analytic with respect to the $x$-variable on $\partial I \times [0,\infty)$~\cite{KukavicaMasmoudiVicolWong14}.
\end{itemize}

However, among the above five settings where the Prandtl equations are known to be locally well-posed, the inviscid limit is known to hold only in the real-analytic setting (b).
This result was established by Sammartino and Caflisch in~\cite{SammartinoCaflisch98b}; see also~\cite{Maekawa14} for a more recent result on vanishing viscosity limit in the analytic setting. In particular, up to our knowledge it is not known whether the inviscid limit \eqref{eq:INVISCID:LIMIT} holds in the Oleinik setting (a), where the solutions have a finite degree of smoothness.

In this paper we prove that the combination of the Oleinik-type condition of no back-flow in the trace of the Euler flow and of a lower bound for the Navier-Stokes vorticity in a boundary layer, imply that the inviscid limit holds.  

A direct connection between the inviscid limit and the one sided-conditions $U\geq 0$ and $\omega|_{\partial \HH} \geq 0$ is provided by the following observation.

\begin{theorem}%[\bf No back-flow \& positive vorticity on the boundary]
Fix $T>0$ and $s>2$, and consider classical solutions $\uns,\ue \in L^\infty(0,T;H^s)$ of \eqref{eq:NS} respectively \eqref{eq:E} with respective boundary conditions \eqref{EQ22} and \eqref{EQ23}. Assume that the trace of the Euler tangential velocity obeys $U(x_1,t) \geq 0$, and that for all $\nu>0$ sufficiently small the trace of the Navier-Stokes vorticity obeys $\omega|_{\partial \HH}\geq0$, for all $x_1\in\RR$ and $t\in [0,T]$.
%Then the inviscid limit \eqref{eq:INVISCID:LIMIT} holds as $\nu \to 0$.
Then 
 \begin{align}
 \| \uns - \ue \|_{L^\infty(0,T;L^2(\HH))}^2 \to 0 
 \label{EQ001}
 \end{align} 
holds as $\nu \to 0$.
\label{thm:boundary}
\end{theorem}
\begin{remark}
\label{rem:Jim}
If follows from the proof of the theorem that instead of assuming $\omega|_{\partial \HH} \geq 0$, we may assume the much weaker condition 
\begin{align}
\omega|_{\partial \HH}  = -\partial_2 \uns_1|_{\partial \HH} \geq - \frac{M_{\nu}(t)}{\nu}
\label{eq:bdry:vort}
\end{align}
for some positive function $M_{\nu}$ which obeys $\int_0^T M_{\nu}(t) dt  \to 0$ as $\nu \to 0$, and  obtain that \eqref{EQ001} holds~\cite{Kelliher14}.
\end{remark}

%Our main result of this paper, Theorem~\ref{thm:layer} below, improves on the Oleinik-type conditions in two aspects. First, it shows that the vorticity does not need to be positive in the boundary layer, it can be in fact almost as negative as $-1/\nu$. Second, the thickness of the boundary layer where the lower bound on $\omega$ is assumed is not $\sqrt{\nu}$, but instead it can be almost as thin as $\nu$. 
%The size of this smaller boundary layer is related to the results of Kato~\cite{Kato84b}, which were later extended by Temam and Wang~\cite{TemamWang97b}. As opposed to these works however, our conditions are one-sided, which is in the spirit of Oleinik's assumptions. 

Our main result of this paper, Theorem~\ref{thm:layer} below, shows that if in a boundary layer almost as thin as $\nu$ the vorticity is not too negative, then the inviscid limit holds. The size of this boundary layer is related to the results of Kato~\cite{Kato84b}, which were later extended by Temam and Wang~\cite{TemamWang97b}. Note however that our conditions are one-sided, which is in the spirit of Oleinik's assumptions.

\begin{theorem}%[\bf No back-flow \& vorticity bounded from below in boundary layer]
Fix $T>0$, $s>2$, and consider classical solutions $\uns,\ue \in L^\infty(0,T;H^s)$ of \eqref{eq:NS} and \eqref{eq:E} respectively with respective boundary conditions \eqref{EQ22} and \eqref{EQ23}. Let $\tau(t) = \min\{t,1\}$ and let $M_{\nu}$ be a positive function which obeys
 \begin{align}
 \int_0^T M_{\nu}(t) dt \to 0 \quad \mbox{as} \quad \nu \to 0. \label{eq:cond:2}
 \end{align}
Define the boundary layer $\Gamma_\nu$ by
\begin{align}
\Gamma_{\nu}(t) = \left\{ (x_1,x_2) \in \HH \colon 0<x_2 \leq \frac{\nu \tau(t)}{C} \log\left(\frac{C}{M_{\nu}(t) \tau(t)}\right)  \right\}  
\label{eq:BL:def}
\end{align} 
where $C= C(\|\ue\|_{L^\infty(0,T;H^s)})>0$ is a sufficiently large fixed positive constant.  
Assume that there is no back-flow in the trace of the Euler tangential velocity, i.e.,
 \begin{align}
 U(x_1,t) \geq 0
 \label{eq:cond:1}
 \end{align} 
for all $x_1 \in \RR$ and $t \in [0,T]$, 
and that for all $\nu$ sufficiently small the ``very negative part'' of the Navier-Stokes vorticity obeys 
\begin{align}
\nu^{(r-1)/r}  \left\| \left( \omega(x_1,x_2,t)+ \frac{M_\nu(t)}{\nu} \right)_- \right\|_{L^r(\Gamma_\nu(t))} \leq  \tau(t)^{1/r} M_{\nu}(t)   \label{eq:cond}
\end{align}
for some $1\leq r \leq \infty$ and all $t\in [0,T]$, where $f_{-} = \min\{f,0\}$. Then the inviscid limit \eqref{eq:INVISCID:LIMIT} holds, with 
the rate of convergence 
\[
\|u - \bar u\|_{L^\infty(0,T;L^2)}^2 = \OO\left(\nu T + \int_0^T M_{\nu}(t) dt\right)
\] 
as $\nu \to 0$.
\label{thm:layer}
\end{theorem}

Note that the above result may be viewed as a one-sided Kato criterion.

\begin{remark}
Since on $\partial \HH$ we have that $\partial_1 \uns_2 =0$, the condition \eqref{eq:cond} on $\omega$ can be replaced by the same condition with $\omega$ replaced by $-\partial_2 \uns_1$. 
\end{remark}

\begin{example}
The shear flow solution $(e^{t \nu \partial_{yy}} v(y),0)$, with $v(0)=0$, and $v'(y)\leq 0$ for $0 \leq y \leq1$ obeys the conditions of Theorem~\ref{thm:layer}.
\end{example}

\begin{remark}
The condition $U \geq 0$ can be ensured for an $\OO(1)$ amount of time if the initial data obeys e.g.~$U_0 \geq \sigma > 0$. However it is not clear that if assuming the initial vorticity obeys $\omega_0 \geq \sigma>0$ implies that \eqref{eq:cond} holds for an $\OO(1)$ time.
\end{remark}

\begin{remark}
We note that Theorem~\ref{thm:layer} also holds in the case of a a bounded domain $\Omega$ with smooth boundary, cf.~Theorem~\ref{thm:boundary:curved_domain} below. The only difference between the inviscid limit on $\HH$ and that on $\Omega$ is that for the later case we need to choose a compactly supported  boundary layer corrector. This is achieved using the argument of~\cite{TemamWang97b}. We refer to Section~\ref{sec:curved} below for details.
\end{remark}

The paper is organized as follows. In Section~\ref{sec:2} we give the proof of Theorem~\ref{thm:boundary}, while in Section~\ref{sec:3} we give the proof of Theorem~\ref{thm:layer}. Lastly, in Section~\ref{sec:curved} we give the main ideas for the proof of Theorem~\ref{thm:boundary:curved_domain}, our main result in the case of a smooth bounded domain.

\section{Proof of Theorem~\ref{thm:boundary}}
\label{sec:2}

Let $v = \uns - \ue$, and $q = \pns-\pe$ be the velocity and the pressure differences respectively. Then $v=(v_1,v_2)$ and $q$ obey the equation
 \begin{align*}
 \partial_t v - \nu \Delta \uns + v \cdot \nabla \ue + \uns \cdot \nabla v + \nabla q = 0
 \end{align*}
with the boundary conditions
 \begin{align*}
 &v_1|_{\partial_\HH} = - \uu, \qquad
 v_2|_{\partial_\HH} = 0.
 \end{align*}
and the initial condition
 \begin{align*}
 v|_{t=0} = 0.
 \end{align*}
The energy identity for the velocity difference then reads
 \begin{align*}
 \frac 12 \frac{d}{dt} \|v\|_{L^2}^2 + \nu \|\nabla \uns \|_{L^2}^2 = - \nu \int_{\HH} \Delta \uns \cdot \ue - \int_{\HH} v \cdot \nabla\ue \; v
 \end{align*}
Using the conditions of the theorem and the given boundary conditions, we get
 \begin{align*}
 - \nu \int_{\HH} \Delta \uns \cdot \ue 
 &=  \nu \int_{\HH} \nabla \uns \cdot \nabla \ue  + \nu \int_{\partial \HH} \partial_2 \uns_1 \; \ue_1 \notag\\
 &= \nu \int_{\HH} \nabla \uns \cdot \nabla \ue - \nu \int_{\partial \HH} \omega \; \uu dx \notag\\
 &\leq  \nu \int_{\HH} \nabla \uns \cdot \nabla \ue \notag\\
 &\leq \nu \|\nabla \uns\|_{L^2}^2 + \frac{\nu}{4} \|\nabla \ue \|_{L^2}^2.
 \end{align*}
We thus obtain
 \begin{align*}
 \frac 12 \frac{d}{dt} \|v\|_{L^2}^2 \leq  \frac{\nu}{4} \|\nabla \ue \|_{L^2}^2 + \|\nabla\ue\|_{L^\infty} \|v\|_{L^2}^2 \leq C \nu + C \|v\|_{L^2}^2
 \end{align*}
 where $C$ is a constant that is allowed to depend on $T$ and $\|\ue\|_{L^\infty(0,T;H^s)}$.
 Recalling that $v(0)=0$, we obtain from the Gr\"onwall Lemma
 \begin{align*}
% \|v(t)\|_{L^2}^2 \leq C \nu t + \frac{\nu}{C} (e^{Ct} - 1 - Ct) \leq C \nu t
 \|v(t)\|_{L^2}^2 \leq C \nu t + C\nu e^{C t} \leq C \nu t
 \end{align*}
which completes the proof. Note that the rate of convergence is $\OO(\nu)$ as $\nu \to 0$.

In order to see that Remark~\ref{rem:Jim} holds, note that under the condition~\eqref{eq:bdry:vort} on the boundary vorticity, one may estimate
\begin{align*}
- \nu \int_{\partial\HH} \omega \uu dx \leq M_{\nu}(t) \int_{\partial \HH} \uu dx  \leq C M_{\nu}(t)
\end{align*}
where in the last inequality we have used a trace inequality. Note moreover that the that the new rate of convergence is ${\mathcal O}(\nu + \int_0^T M_{\nu}(t) dt)$.

\section{Proof of Theorem~\ref{thm:layer}}
\label{sec:3}

In the spirit of~\cite{Kato84b}, the proof is based on constructing a suitable boundary layer corrector $\phi$ to account for the mismatch between the Euler and Navier-Stokes boundary conditions.
Note however that the Kato's corrector $\phi=(\phi_1,\phi_2)$ is not suitable here due to the change of sign of
$\phi_1$.

\subsection*{The boundary layer corrector}

We fix $\psi \colon [0,\infty) \to [0,\infty)$ to be a $C_0^\infty$ function approximating $\chi_{[1,2]}$, supported in $[1/2,4]$, which is non-negative and has mass $\int \psi(z) dz = 1$. Recall that $\tau(t) = \min\{t , 1\}$.

For $\alpha \in (0,1]$, to be chosen later, we introduce
 \begin{align*} 
 \phi(x_1,x_2,t) = (\phi_1(x_1,x_2,t),\phi_2(x_1,x_2,t))
 \end{align*}
where
 \begin{align}
   \phi_1(x_1,x_2,t) &=  - \uu(x_1,t)\left(e^{-\fractext{x_2}{\alpha \tau(t)}} - \alpha \tau(t) \psi(x_2) \right) \label{eq:phi:1} \\
   \phi_2(x_1,x_2,t) &=   \alpha \tau(t) \partial_{1} \uu(x_1,t) \left( \left( 1 - \int_0^{x_2} \psi(y) dy \right) - e^{-\fractext{x_2}{\alpha \tau(t)}} \right) \label{eq:phi:2}
 \end{align}
and 
 \begin{align*}
 \phi(x_1,x_2,0) = \phi_0(x_1,x_2)= 0.
 \end{align*}
Observe that we have 
$\phi_1\to0$ as $x_2\to\infty$ exponentially, and
 \begin{align*}
 &\phi_1(x_1,0,t) = - \uu(x_1,t) 
 \nonumber\\&
 \phi_2(x_1,0,t) = 0.
 \end{align*}
In particular, note that
 \begin{align*} 
 \ue + \phi = 0 \quad \mbox{on}\quad \partial \HH.
 \end{align*}
Equally importantly, the corrector is divergence free
 \begin{align*} 
 \nabla \cdot \phi = 0
 \end{align*}
which allows us not to deal with the pressure when performing energy estimates.

Throughout the proof, we shall also use the bounds
 \begin{align*} 
  \| \phi_1 \|_{L^p} &\leq C (\alpha \tau)^{1/p} + C \alpha \tau \leq C (\alpha \tau)^{1/p}
 \end{align*}
and
 \begin{align*}
  &
  \| \partial_1 \phi_1 \|_{L^p} \leq C (\alpha \tau)^{1/p} 
  \nonumber\\&
  \| \partial_2 \phi_1 \|_{L^p} \leq C (\alpha \tau)^{1/p - 1}
 \end{align*}
for any $1\leq p \leq \infty$, with
 \begin{align*} 
 &\|\phi_2 \|_{L^p} \leq C \alpha \tau ( 1 +(\alpha \tau)^{1/p} ) \leq C \alpha \tau 
 \nonumber\\&
 \|\partial_1 \phi_2 \|_{L^p} \leq C \alpha \tau (1 + (\alpha \tau)^{1/p} )\leq C \alpha \tau
 \end{align*}
since $\alpha \tau \leq \alpha \leq 1$.
Here and throughout the proof, the constant $C$ is allowed to depend on various norms of $\uu$ and $\ue$ (which we do not keep track of), but not on norms of $\uns$.

\subsection*{Energy equation} As before, define the velocity and pressure differences by 
 \begin{align*} 
 &v= \uns - \ue
 \nonumber\\&
 q = \pns - \pe.
 \end{align*} 
Subtracting \eqref{eq:E} from \eqref{eq:NS} we arrive at
 \begin{align} 
 \partial_t (v-\phi) - \nu \Delta \uns + v \cdot \nabla \ue + \uns \cdot \nabla v + \nabla q  + \partial_t \phi = 0.
 \label{eq:v-phi}
 \end{align}
Since $v-\phi=u-\ue-\phi=0$ on $\partial \HH$, we may multiply \eqref{eq:v-phi} by $v-\phi$ and integrate by parts to obtain
 \begin{align} 
 &\frac 12 \frac{d}{dt} \| v- \phi\|_{L^2}^2 + \nu \|\nabla u\|_{L^2}^2
 \nonumber\\&\qquad
 = \nu \int \nabla \uns \nabla  \phi 
    - \int \uns \nabla \phi \uns \notag
    \nonumber\\&\qquad\qquad
    + \Bigg( \nu \int \nabla \uns \nabla \ue  - \int (v-\phi) \nabla \ue (v-\phi) \notag
%    \\ &\qquad 
    - \int \phi \nabla \ue (v-\phi) + \int \uns \nabla \phi \ue  - \int \partial_t \phi (v-\phi) \Bigg)\notag
    \\&
     \qquad = I_1 + I_2 + R
 \label{eq:ODE:1}
 \end{align}
where we used
$-\int v\nabla \ue(v-\phi)=-\int(v-\phi)\nabla\ue(v-\phi)-\int\phi\nabla\ue(v-\phi)$
and
$-\int u\nabla v(v-\phi)=-\int u\nabla\phi(v-\phi)=-\int u\nabla\phi(u-\ue)$.
The terms $I_1$ and $I_2$ give the main contributions, while the $R$ term is in some sense a remainder term. The assumptions on the sign of $\uu$ and on the very negative part of $\omega$ come into play when bounding $I_1$.

\subsection*{Estimate for $I_1$}

We decompose $I_1$ as 
 \begin{align*} 
 I_1 = \nu \int \nabla \uns \nabla \phi 
 = \nu \int \partial_2 \uns_1 \partial_2 \phi_1 + \sum_{(i,j)\neq (1,2)} \nu \int \partial_i \uns_j \partial_i \phi_j = I_{11} + I_{12}.
 \end{align*}
In order to estimate $I_{12}$, we consider the three possible combinations of $(i,j)\neq (1,2)$. We have
 \begin{align*} 
 \nu \int \partial_1 \uns_1 \partial_1 \phi_1 &\leq \frac{\nu}{4} \|\partial_1 \uns_1\|_{L^2}^2 +  \nu \|\partial_1 \phi_1\|_{L^2}^2 \leq \frac{\nu}{4} \|\partial_1 \uns_1\|_{L^2}^2 + C \nu (\alpha \tau) \\
 \nu \int \partial_1 \uns_2 \partial_1 \phi_2 &\leq \frac{\nu}{4} \|\partial_1  \uns_2\|_{L^2}^2 + \nu \|\partial_1 \phi_2 \|_{L^2}^2 \leq \frac{\nu}{4} \|\partial_1  \uns_2\|_{L^2}^2 + C \nu (\alpha \tau)^2 \\
 \nu \int \partial_2 \uns_2 \partial_2 \phi_2 &\leq \frac{\nu}{4} \|\partial_2 \uns_2\|_{L^2}^2 + \nu \|\partial_2 \phi_2\|_{L^2}^2 \leq \frac{\nu}{4} \|\partial_2 \uns_2\|_{L^2}^2 + C \nu (\alpha \tau)
 \end{align*}
for some sufficiently large $C$, which shows that 
 \begin{align}
 I_{12} \leq \frac{\nu}{4} \|\nabla u\|_{L^2}^2 + C \nu (\alpha \tau).
 \label{eq:I12}
 \end{align}

The main contribution to $I_1$ comes from the term $I_{11}$, which we bound next. Let $\beta$ be the thickness of the boundary layer where the assumption on the very negative part of $\omega = \partial_1 \uns_2 - \partial_2 \uns_1$ is imposed. That is, for some $\beta \in (\alpha, 1/4]$ and $M>0$, to be specified below, we use the bound
 \begin{align}
 \omega(x_1,x_2,t) \geq - \frac{M}{\nu} + \tilde\omega(x_1,x_2,t)
   \comma  (x_1,x_2) \in \Gamma_{\beta} =  \RR \times (0,\beta)
   , \quad t \in [0,T],
 \label{eq:omega:LWB}
 \end{align}
where we have denoted
\begin{align*}
\tilde\omega(x_1,x_2,t) = \min\left \{ \omega(x_1,x_2,t) + \frac{M}{\nu}, 0 \right\} \leq 0.
\end{align*}
Next, we decompose 
 \begin{align*}
 I_{11} 
 &= \nu \int_{\HH} \partial_2 \uns_1 \partial_2 \phi_1 
 = - \nu \int_{\Gamma_\beta} \omega \partial_2 \phi_1 - \nu \int_{\Gamma_\beta^C} \omega \partial_2 \phi_1 + \nu \int_{\HH} \partial_1 \uns_2 \partial_2 \phi_1 \notag\\
 &= I_{111} + I_{112} + I_{113}.
 \end{align*}
The assumptions \eqref{eq:cond:1} and \eqref{eq:omega:LWB} are only be used to estimate $I_{111}$.
By construction of the corrector in \eqref{eq:phi:1}--\eqref{eq:phi:2}, we have the explicit formula
 \begin{align} 
 \partial_2 \phi_1(x_1,x_2,t) =  \frac{1}{\alpha \tau} \uu(x_1,t) e^{-\fractext{x_2}{\alpha \tau}} - \alpha \tau \uu(x_1,t) \psi'(x_2)  
 \label{eq:corrector:D2}
 \end{align}
for all $(x_1,x_2) \in \HH$ and $t\in [0,T]$. 
In view of the no-back flow condition $\uu \geq 0$ of \eqref{eq:cond:1} and the bound \eqref{eq:omega:LWB} on $\omega$ in $\Gamma_\beta$, for any $r \in [1,\infty]$ we have the estimate
 \begin{align}
 I_{111} 
 &= - \frac{\nu}{\alpha \tau} \int_{\Gamma_\beta} \omega(x_1,x_2,t) \uu(x_1,t) e^{-\fractext{x_2}{\alpha \tau}} dx_1 dx_2 \notag\\
 &\qquad + \nu \alpha \tau \int_{\Gamma_\beta} \omega(x_1,x_2,t) \uu(x_1,t)  \psi'(x_2) dx_1 dx_2 \notag\\
 &\leq \frac{M }{\alpha \tau} \int_{x_2<\beta} \uu(x_1,t) e^{-\fractext{x_2}{\alpha \tau}} dx_1 dx_2 + \frac{\nu}{\alpha \tau} \int_{x_2<\beta} (- \tilde\omega(x_1,x_2,t)) \uu(x_1,t) e^{-\fractext{x_2}{\alpha \tau}} dx_1 dx_2\notag\\
 &\qquad + \nu \alpha \tau \int_{x_2<\beta} \omega(x_1,x_2,t) \uu(x_1,t)  \psi'(x_2) dx_1 dx_2 \notag\\
 &\leq  C M (1 - e^{-\beta/\alpha \tau})+ \frac{C \nu}{\alpha \tau} \|\tilde\omega\|_{L^r(\Gamma_\beta)} \left(  \frac{(r-1)\alpha \tau}{r}  (1 - e^{-r\beta/(r-1)\alpha \tau}) \right)^{(r-1)/r} + C \nu \alpha \tau \|\nabla \uns\|_{L^2} \notag\\
 &\leq \frac{\nu}{12} \|\nabla \uns\|_{L^2}^2 + C M  + C \nu (\alpha \tau)^{-1/r}  \|\tilde\omega\|_{L^r(\Gamma_\beta)} + C \nu (\alpha \tau)^2 
 \label{eq:I111}
 \end{align}
for a sufficiently large $C$. For the outer layer term $I_{112}$ we have
 \begin{align}
 I_{112} &= - \nu \int_{x_2>\beta} \omega(x_1,x_2,t) \uu(x_1,t) \left( \frac{1}{\alpha \tau} e^{-\fractext{x_2}{\alpha \tau}} - \alpha \tau \psi'(x_2) \right) dx_1 dx_2 \notag\\
 &\leq C \frac{\nu e^{-\fractext{\beta}{\alpha \tau}}}{(\alpha \tau)^{1/2}}  \|\nabla u\|_{L^2} + \nu \alpha \tau \|\nabla u\|_{L^2} \notag\\
 &\leq \frac{\nu}{12} \| \nabla u\|_{L^2}^2 + C \frac{\nu}{\alpha \tau} e^{-\fractext{2\beta}{\alpha \tau}} +  C \nu (\alpha \tau)^2. 
 \label{eq:I112}
 \end{align}
Lastly, for $I_{113}$ we integrate by parts once in $x_2$ by using that $\partial_1 u_2 = 0 $ on $\partial \HH$, use that $\nabla \cdot \uns = 0$ on ${\mathbb H}$, and then integrate by parts in $x_1$ to obtain
 \begin{align}
 I_{113} &= \nu \int \partial_1 \uns_2 \partial_2 \phi_1  = - \nu \int \partial_{12} \uns_2 \phi_1 = \nu \int \partial_{11} \uns_1 \phi_1 \notag\\
 &= - \nu \int \partial_1 \uns_1 \partial_1 \phi_1 \leq \frac{\nu}{12} \|\nabla u\|_{L^2}^2 + C \nu (\alpha \tau).
 \label{eq:I113}
 \end{align}
Combining \eqref{eq:I111}, \eqref{eq:I112}, and \eqref{eq:I113}, we arrive at 
 \begin{align} 
 I_{11} \leq  \frac{\nu}{4} \| \nabla u \|_{L^2}^2 + C M  + C \nu (\alpha \tau)^{-1/r}  \|\tilde\omega\|_{L^r(\Gamma_\beta)} + C \frac{\nu}{\alpha \tau} e^{-\fractext{2\beta}{\alpha \tau}} +  C \nu (\alpha \tau).
 \label{eq:I11}
 \end{align}
We summarize \eqref{eq:I12} and \eqref{eq:I11} as 
 \begin{align}
 I_1 \leq  \frac{\nu}{2} \| \nabla u \|_{L^2}^2 + C M  + C \nu (\alpha \tau)^{-1/r}  \|\tilde\omega\|_{L^r(\Gamma_\beta)} + C \frac{\nu}{\alpha \tau} e^{-\fractext{2\beta}{\alpha \tau}} +  C \nu (\alpha \tau) \label{eq:I1}
 \end{align}
for a suitable constant $C$ which may depend on norms of $\uu$ and $\ue$.

\subsection*{Estimate for $I_2$}

In order to treat $I_2$, we 
use $\partial_{2}\phi_2=-\partial_{1}\phi_1$ and decompose
 \begin{align*} 
 I_2 = - \int \uns_2 \partial_2 \phi_1 \uns_1 - \int \uns_1 \partial_1 \phi_2 \uns_2 +  \int (\uns_2^2 - \uns_1^2) \partial_1 \phi_1 = I_{21} + I_{22} + I_{23}.
 \end{align*}
We note that upon integration, we have that for any $j \in \{1,2\}$ 
 \begin{align*} 
 \int_0^\infty \uns_j(x_1,x_2)^2  e^{-\fractext{x_2}{\alpha \tau}} dx_2 
 &= 2 \alpha  \tau \int_0^\infty \uns_j(x_1,x_2) \partial_2 \uns_j(x_1,x_2) e^{-\fractext{x_2}{\alpha \tau}}   dx_2   \notag\\
 &\leq C \alpha \tau \left( \int_0^\infty \uns_j(x_1,x_2)^2  e^{-\fractext{x_2}{\alpha \tau}}   dx_2 \right)^{1/2}  \|\partial_2 u_j e^{-x_2/\alpha\tau}\|_{L^2}
 \end{align*}
and thus
 \begin{align*} 
 \|u_j e^{-\fractext{x_2}{2 \alpha \tau}}  \|_{L^2} \leq C \alpha \tau \|\partial_2 u_j\|_{L^2}.
 \end{align*}
Using the above estimate we obtain 
 \begin{align*} 
 I_{21} &= - \frac{1}{\alpha \tau} \int \uns_2 \uu(x_1,t)   e^{-\fractext{x_2}{\alpha \tau}}  \uns_1 +  \alpha\tau \int \uns_2 \uu(x_1,t) \psi'(x_2)  \uns_1 \notag \\
 &\leq
   \frac{C}{\alpha\tau}
   \Vert u_1 e^{-x_2/2\alpha\tau}\Vert_{L^2}
   \Vert u_2 e^{-x_2/2\alpha\tau}\Vert_{L^2}
   +
   C\alpha\tau
   \Vert u_1\Vert_{L^2}
   \Vert u_2\Vert_{L^2}
 \\
 &\leq C \alpha \tau \| \nabla \uns\|_{L^2}^2 + C \alpha \tau \| \uns \|_{L^2}^2 \notag \\
 &\leq C \alpha \tau \| \nabla \uns\|_{L^2}^2 + C \alpha  \tau
 \end{align*}
due to the energy inequality 
 \begin{align*}
 \| \uns\|_{L^2} \leq \| \uns_0\|_{L^2} =  \| \ue_0\|_{L^2} \leq C
 \end{align*} 
which is a viscosity-independent bound.
Similarly, 
 \begin{align*} 
 I_{22} \leq  \|\uns\|_{L^2}^2 \|\partial_1 \phi_2\|_{L^\infty} \leq C \alpha \tau
 \end{align*}
and
 \begin{align*}
 I_{23} \leq C (\alpha \tau)^2 \| \nabla \uns\|_{L^2}^2 + C \alpha \tau. 
 \end{align*}
Using that $\alpha \tau \leq 1$, we arrive at 
 \begin{align} 
 I_2 \leq  C \alpha \tau \| \nabla \uns\|_{L^2}^2 + C \alpha  \tau
 \label{eq:I2}
 \end{align}
where as usual, the constant $C$ is allowed to depend on  the Euler flow.

\subsection*{Estimate for the remainder terms $R$}
Using
%\cor
  \begin{equation*}
    \int \uns\nabla\phi\ue
    =-\int \uns\nabla\ue\phi
    =-\int(v-\phi)\nabla\ue\phi-\int\phi\nabla\ue\phi-\int\ue\nabla\ue\phi, 
  \end{equation*}
%\cob
we may rewrite the remainder term as
 \begin{align*} 
 R &= \nu \int \nabla \uns \nabla \ue  - \int (v-\phi) \nabla \ue (v-\phi) - \int \phi \nabla \ue (v-\phi)  \notag\\
 & \qquad - \int (v-\phi) \nabla \ue \phi  - \int \phi \nabla \ue \phi -\int \ue\nabla\ue\phi- \int \partial_t \phi (v-\phi) \notag\\
 &=: R_1 + R_2 + R_3 + R_4 + R_5 + R_6 + R_7.
 \end{align*}
Using the available bounds on the corrector $\phi$, we have the bounds
 \begin{align*}
 R_1 &\leq \nu \|\nabla \uns \|_{L^2} \|\nabla \ue\|_{L^2} \leq \frac{\nu}{4} \|\nabla \uns\|_{L^2}^2 + C \nu \notag\\
 R_2 &\leq \|v-\phi\|_{L^2}^2 \|\nabla \ue\|_{L^\infty} \leq C \|v-\phi\|_{L^2}^2 \notag\\
 R_3 &\leq \|\phi\|_{L^2} \|\nabla \ue\|_{L^\infty} \|v-\phi\|_{L^2} \leq C (\alpha \tau)^{1/2} \|v-\phi\|_{L^2} \leq C \|v-\phi\|_{L^2}^2 + C \alpha \tau\notag\\
 R_4 &\leq \|v-\phi\|_{L^2} \|\nabla \ue\|_{L^\infty} \|\phi\|_{L^2} \leq C \|v-\phi\|_{L^2}^2 + C \alpha \tau\notag\\
 R_5 &\leq \|\phi\|_{L^2}^2 \|\nabla \ue\|_{L^\infty} \leq C \alpha \tau \notag\\
 R_6 &\leq C\|\phi\|_{L^1}\le C\alpha\tau\\
 R_7 &\leq \|\partial_t \phi\|_{L^2} \|v-\phi\|_{L^2} \leq C ( (\alpha \tau)^{1/2} + \alpha ) \|v-\phi\|_{L^2} \leq C \|v-\phi\|_{L^2}^2 + C ( (\alpha \tau) + \alpha^2)
 \end{align*}
which may be summarized as
 \begin{align} 
 R \leq \frac{\nu}{4} \|\nabla \uns\|_{L^2}^2 + C \|v-\phi\|_{L^2}^2 + C \nu + C \alpha^2 + C \alpha \tau
 \label{eq:R}
 \end{align}
and again, $C$ depends on various norms of $\uu$ and  $\ue$.

\subsection*{Conclusion of the proof}
Combining \eqref{eq:ODE:1}, with \eqref{eq:I1}, \eqref{eq:I2}, and \eqref{eq:R}, we thus arrive at the bound
 \begin{align}
 & \frac 12 \frac{d}{dt} \| v- \phi\|_{L^2}^2 + \left(\frac{\nu}{4}- C \alpha \tau \right) \|\nabla \uns\|_{L^2}^2 \notag\\
 &\qquad  \leq  C \|v-\phi\|_{L^2}^2 + C \nu + C \alpha^2 + C \alpha \tau  + C M  + C \frac{\nu}{\alpha \tau} e^{- \fractext{2\beta}{\alpha \tau}} + C \nu (\alpha \tau)^{-1/r}  \|\tilde\omega\|_{L^r(\Gamma_\beta)}
 \label{eq:ODE:2}
 \end{align}
for all $\nu \leq 1$.
To conclude the proof, we first choose
 \begin{align} 
 \alpha = \frac{\nu}{C}
 \label{eq:alpha:def}
 \end{align}
for a sufficiently large $C$, since $\tau(t) \leq 1$, we obtain from \eqref{eq:ODE:2} that
 \begin{align}
 & \frac 12 \frac{d}{dt} \| v- \phi\|_{L^2}^2  \leq  C \|v-\phi\|_{L^2}^2 + C \nu  + C M + C \tau^{-1} e^{- \fractext{2C\beta}{\nu \tau}} + C \nu^{(r-1)/r}  \tau^{-1/r}  \|\tilde\omega\|_{L^r(\Gamma_\beta)}
 \label{eq:ODE:3}
 \end{align}
for $\nu \leq 1$. Next, per our assumption  \eqref{eq:BL:def}, define the boundary layer thickness by 
 \begin{align*}
% \beta = \frac{\alpha \tau}{2} \log \left( \frac{1}{\alpha \tau} \right)
 \beta = \frac{\nu \tau}{2C} \log \left( \frac{1}{M \tau} \right)
 \end{align*}
where $\alpha$ is given by \eqref{eq:alpha:def}, and obtain from assumption \eqref{eq:ODE:3} that
 \begin{align}
 & \frac 12 \frac{d}{dt} \| v- \phi\|_{L^2}^2  \leq  C \|v-\phi\|_{L^2}^2 + C \nu  + C M  
 \label{eq:ODE:4}
 \end{align}
Since 
 \begin{align*} 
 \|v_0-\phi_0\|_{L^2}^2 = \|\phi_0\|_{L^2}^2 = 0 
 \end{align*}
the Gr\"onwall inequality applied to \eqref{eq:ODE:4} yields
 \begin{align}
 \|\uns(t) - \ue(t) \|_{L^2}^2 
 &\leq \|\phi(t)\|_{L^2}^2 + \|v(t) - \phi(t)\|_{L^2}^2 \notag\\
 &\leq C \nu \tau(t) + C T e^{CT} \left( \nu t + \int_0^t M(s) ds \right)
 \label{eq:ODE}
 \end{align}
for all $t \in [0,T]$. The assumption \eqref{eq:cond:2} yields that the right side of \eqref{eq:ODE} converges to $0$ as $\nu \to 0$, uniformly in $t \in [0,T]$, which concludes the proof.

\section{Curved domains}
\label{sec:curved}
In this section, we show that
Theorem~\ref{thm:layer}
holds also in the case of a bounded domain
$\Omega$  with a smooth boundary and with an outer normal $n$.

\begin{theorem}%[\bf No back-flow \& vorticity bounded from below in boundary layer, curved domain]
Fix $T>0$ and $s>2$, and consider classical solutions $\uns,\ue \in
L^\infty(0,T;H^s)$ of \eqref{eq:NS} respectively \eqref{eq:E} in $\Omega$ with respective boundary conditions $u|_{\partial\Omega}=0$ and $\ue\cdot n|_{\partial\Omega}=0$. Let $M_{\nu}$ be a positive function such that
\begin{align}
 \int_0^T M_{\nu}(t) dt \to 0 \quad \mbox{as} \quad \nu \to 0, \label{eq:cond:2:curved}
\end{align}
and define the boundary layer
 \begin{align}
 \Gamma_{\nu}(t) =  \left\{
   x\in\Omega:
%      0<x_2 \leq \frac{\nu \tau(t)}{C} \log\left(\frac{C}{\nu \tau(t)}\right)
      0<\dist(x,\partial\Omega) \leq \frac{\nu \min\{t,1\}}{C} \log\left(\frac{C}{M_{\nu}(t) \min\{t,1\}}\right)
  \right\}  
 \label{eq:BL:def:curved}
 \end{align} 
 where $C= C(\|\ue\|_{L^\infty(0,T;H^s)})>0$ is a sufficiently large fixed constant. 
Assume that the trace of the Euler tangential velocity (cf.~\eqref{EQ02} below) is nonnegative and that for all $\nu>0$ sufficiently small the Navier-Stokes vorticity obeys 
 \begin{align}
\nu^{(r-1)/r} \left\| \min\left\{ \omega(\cdot,t) + \frac{M_\nu(t)}{\nu}, 0 \right\} \right\|_{L^r(\Gamma_\nu(t))} \leq M_\nu(t) \min\{ t,1\}^{1/r} \label{eq:cond:curved}
 \end{align}
for some $1\leq r \leq \infty$ and all $t \in [0,T]$. Then the inviscid limit \eqref{eq:INVISCID:LIMIT} holds as $\nu \to 0$, with
a rate of convergence 
proportional to $\nu T +\int_0^T M_{\nu}(t) dt$. 
\label{thm:boundary:curved_domain}
\end{theorem}

{\begin{proof}
We follow the notation and ideas from
\cite{TemamWang97b}. Note however that we need
to modify the corrector 
since the one in \eqref{eq:phi:1}--\eqref{eq:phi:2}
is not compactly supported.

As in \cite{TemamWang97b}, let
$(\xi_1,\xi_2)$ denote the
orthogonal coordinate system defined in a sufficiently small
neighborhood of the boundary
  \begin{equation}
   \bigl\{
    x=(x_1,x_2)\in{\mathbb R}^2:
    \dist(x,\partial\Omega)\le\delta
   \bigr\}
   \label{EQ00}
  \end{equation}
where $\delta>0$.
%\cor with the positive orientation\cob. 
Here
$\xi_2$ denotes the distance to the boundary $\partial\Omega$.
For simplicity of notation, we assume that
$\partial\Omega$ consists of the connected
smooth Jordan curve---the modification to the general case
of finitely many Jordan curves can be done similarly.
Then we may use the notation from
\cite{TemamWang97b,Batchelor70}; in particular,
  \begin{equation}
   dx_1 dx_2
   =h(\xi_1,\xi_2) d\xi_1^2
   +d\xi_2^2
   .
   \label{EQ01}
  \end{equation}
Denote by $e_1(\xi_1,\xi_2)$ and
          $e_2(\xi_1,\xi_2)$ the local 
basis in the directions of $\xi_1$ and $\xi_2$
respectively. Also, write
  \begin{equation}
   U(\xi_1,t)
   =\ue(\xi_1,0,t)\cdot e_2(\xi_1,0)
   \label{EQ02}
  \end{equation}
for the trace of the Euler flow $\ue$.
Let $\eta(y)\in C_{0}^{\infty}({\mathbb R},[0,1])$ denote
the function which equals~$1$ in a neighborhood
of $(-\infty,-\delta]\cup[\delta,\infty)$, and
let $\psi\in C_{0}^{\infty}({\mathbb R},[0,1])$
be a function supported in
the interval $(\delta/2,\delta)$ such that
$\int \psi=1$. Then define the corrector
  \begin{equation}
   \phi(\xi_1,\xi_2,t)
   =\curl\psi(\xi_1,\xi_2,t)
   \label{EQ03}
  \end{equation}
where
  \begin{align}
   \psi(\xi_1,\xi_2,t)
   = - U(\xi_1,t)
     \int_{0}^{\xi_2}
       \exp\left(
            -\frac{y}{\alpha\tau}
           \right)
           \eta(y)
        \,dy
    + \gamma(t) U(\xi_1,t)
      \int_{0}^{\xi_2}\psi(y)\,dy
    .
   \label{EQ04}
  \end{align}
The parameter $\gamma=\gamma(t)$ is chosen so that
$\psi$ vanishes on $[\delta,\infty)$.
Using $\int\psi=1$, this holds if
  \begin{equation}
   \gamma(t)
   = \int_{0}^{\delta}
            \exp\left(
            -\frac{y}{\alpha\tau}
           \right)
           \eta(y)
        \,dy
    .
   \label{EQ05}
  \end{equation}
Note that $\gamma$ does not depend on $(\xi_1,\xi_2)$
and that we have
  \begin{equation}
   \gamma(t)
   =\alpha\tau(t)
   +{\mathcal O}((\alpha\tau)^{3})
   .   
   \label{EQ06}
  \end{equation}
From \cite{Batchelor70}, recall the formulas
  \begin{equation}
   \ddiv u
   = \frac{1}{h}
     \frac{\partial u_1}{\partial\xi_1}
     + \frac{1}{h}
       \frac{\partial}{\partial\xi_2}(h u_2)
   \label{EQ07}
  \end{equation}
and
  \begin{equation}
   \curl f
   = \frac{\partial f}{\partial\xi_2} e_1
     - \frac{1}{h}
       \frac{\partial}{\partial\xi_1}(h f) e_2
   \label{EQ08}
  \end{equation}
for every vector function $u$ and scalar function $f$
respectively.
Thus we have
  \begin{align}
   \phi_1
   = - U(\xi_1,t)
     \exp\left(
          -\frac{\xi_2}{\alpha\tau}\eta(\xi_2)
         \right)
     + \gamma U(\xi_1,t)\psi(\xi_2)
   \label{EQ10}
  \end{align}
and
  \begin{equation}
   \phi_2
   =
   \frac{1}{h}
   \frac{\partial}{\partial\xi_1}(h U)
     \int_{0}^{\xi_2}
       \exp\left(
            -\frac{y}{\alpha\tau}
           \right)
           \eta(y)
        \,dy
     - \frac{\gamma}{h}
       \frac{\partial}{\partial\xi_1}(h U)
       \int_{0}^{\xi_2}\psi(y)\,dy
    .
   \label{EQ09}
  \end{equation}
As in the previous sections, we have
  \begin{align}
   \frac12
   \frac{d}{dt}
   \Vert v-\phi\Vert_{L^2}^2
   +\nu\Vert\nabla u\Vert_{L^2}^2
   = I_1+I_2+R
   \label{EQ11}
  \end{align}
where
  \begin{equation}
   I_1
   =\nu\int\nabla u \nabla\phi
   \label{EQ12}
  \end{equation}
and
  \begin{equation}
   I_2
   =-\int \uns\nabla \phi \uns
   \label{EQ13}
  \end{equation}
with
  \begin{align}
   R
   &= \nu\int\nabla \uns\nabla\ue
   -\int (v-\phi)\nabla\ue(v-\phi)
   -\int\phi\nabla\ue(v-\phi)
   \nonumber\\&\indeq
   -\int(v-\phi)\nabla\ue\phi
   -\int\phi\nabla\ue\phi
   -\int\ue\nabla\ue\phi
   -\int\partial_{t}\phi(v-\phi)
%   \nonumber\\&
%   =R_1 + R_2 + R_3 + R_4 + R_5 + R_6+ R_7
   .
   \label{EQ14}
  \end{align}
Here, we treat the term 
  \begin{align}
   I_1 = \nu\int\nabla u\nabla\phi
     &
     =\nu
     \int
     \left(
      \frac{1}{h}\frac{\partial\uns_1}{\partial\xi_1}e_1
      +\frac{\partial\uns_1}{\partial\xi_2}e_2
     \right)
     \left(
      \frac{1}{h}
      \frac{\partial\phi_1}{\partial\xi_1}e_1
      +\frac{\partial\phi_1}{\partial\xi_2}e_2
     \right)h\,d\xi_1 d\xi_2
     \nonumber\\&\indeq
     +\nu
     \left(
      \frac{1}{h}\frac{\partial\uns_2}{\partial\xi_1}e_1
      +\frac{\partial\uns_2}{\partial\xi_2}e_2
     \right)
     \left(
      \frac{1}{h}
      \frac{\partial\phi_2}{\partial\xi_1}e_1
      +\frac{\partial\phi_2}{\partial\xi_2}e_2
     \right)h\,d\xi_1 d\xi_2
     .
   \label{EQ15}
  \end{align}
while the rest are estimated
similarly to \cite{TemamWang97b}.
We write the far right side of \eqref{EQ15}
as $I_{11}+I_{12}$  where
  \begin{align}
   I_{11}
   = \nu
     \int \frac{\partial u_1}{\partial\xi_2}
          \frac{\partial \phi_1}{\partial\xi_2}h
   \label{EQ16}
  \end{align}
and $I_{12}$ is the sum of the other three terms.
Then
  \begin{align}
   I_{11}
    &= \nu
       \int \frac{\partial(h \uns_1)}{\partial\xi_2}
            \frac{\partial\phi_1}{\partial\xi_2}
      -\nu
      \int \uns_1 
           \frac{\partial h}{\partial\xi_2}
           \frac{\partial\phi_1}{\partial\xi_2}
    \nonumber\\&
     \qquad -\nu\int_{\xi_2<\beta}
        \omega\frac{\partial\phi_1}{\partial\xi_2}
     -\nu\int_{\xi_2\ge\beta}
        \omega\frac{\partial\phi_1}{\partial\xi_2}
     +\nu
      \int\frac{\partial\uns_2}{\partial\xi_1}
          \frac{\partial\phi_1}{\partial\xi_2}
     -\nu
     \int\uns_1
         \frac{\partial h}{\partial\xi_2}
         \frac{\partial\phi_1}{\partial\xi_2}
   \nonumber\\&
   =
   I_{111}
   +I_{112}
   +I_{113}
   +I_{114}
   \label{EQ17}
  \end{align}
where
  \begin{equation}
   \omega
   = \frac{\partial \uns_2}{\partial\xi_1}
     -\frac{\partial(h \uns_1)}{\partial\xi_2}
   \label{EQ18}
  \end{equation}
is the vorticity.
The terms $I_{111}$ and $I_{112}$ are estimated
the same way as in the flat case. 
For the third term $I_{113}$, we integrate
by parts and obtain
  \begin{equation}
   I_{113}   
   = \nu\int \frac{\partial\uns_2}{\partial\xi_1}
             \frac{\partial\phi_1}{\partial\xi_2}
   = - \nu
       \int\frac{\partial^2\uns_2}{\partial\xi_1\partial\xi_2}\phi_1
   =\nu\int \frac{\partial\uns_2}{\partial\xi_2}
            \frac{\partial\phi_1}{\partial\xi_1}
   .
   \label{EQ19}
  \end{equation}
In the last step we used that
$\partial u_2/\partial\xi_2$ vanishes on the boundary,
which holds since
  \begin{align}
   \frac{\partial\uns_2}{\partial\xi_2}
    =
    \ddiv u
    - \uns_2 \frac{1}{h}\frac{\partial h}{\partial\xi_2}
    - \frac{1}{h}  \frac{\partial\uns_1}{\partial\xi_1}
   .
   \label{EQ20}
  \end{align}
The rest of the terms are treated analogously as in the flat case, and
we obtain
  \begin{equation}
   I_1
   \leq  
    \frac{\nu}{2} \| \nabla u \|_{L^2}^2 
    + C M  + C \nu (\alpha \tau)^{-1/r}  \| \min\{ \omega + M \nu^{-1},0\}\|_{L^r(\Gamma_\beta)}
    + C \frac{\nu}{\alpha \tau} e^{-\fractext{2\beta}{\alpha \tau}} 
    + C \nu (\alpha \tau)   
   .
   \label{EQ21}
  \end{equation}
The terms $I_2$ and $R$ are estimated as in the flat case
(see also \cite{TemamWang97b}), and we thus omit further details.
\end{proof}

\section*{Acknowledgments} 
The authors would like to thank J.P.~Kelliher for suggesting the result in Remark~\ref{rem:Jim}. 
The work of PC was supported in part by the NSF grants DMS-1209394 and DMS-1265132, 
IK was supported in part by the NSF grant DMS-1311943, 
while the work of VV was supported in part by the NSF grant DMS-1211828.

%%%%%%%%%%%%%%%% BIBLIOGRAPHY %%%%%%%%%%%%%%%%%%%%%%%%%%%%%%%%%%%%%%%%%%%

%\bibliographystyle{alpha}
%\bibliography{VladBib}

\end{document}